\documentstyle[aps,amsfonts,epsfig,preprint,pra]{revtex}

\begin{document} 
\title{\bf Extended Bell and Stirling numbers from hypergeometric exponentiation}
\author{J.-M. Sixdeniers,\footnotemark[1] \footnotetext[1]{e-mail : sixdeniers@lptl.jussieu.fr} K.A Penson\footnotemark[5] \footnotetext[5]{e-mail : penson@lptl.jussieu.fr}
and A.I Solomon\footnotemark[6]\footnotemark[4] \footnotetext[4]{e-mail :
a.i.solomon@open.ac.uk}}\footnotetext[6]{Permanent address: Quantum Processes Group, Open University, Milton Keynes, MK7 6AA, United Kingdom.}
\address{Universit\'{e} Pierre et Marie Curie, Laboratoire de Physique Th\'{e}orique des Liquides, Tour 16, $5^{i\grave{e}me}$ \'{e}tage, 4, place Jussieu, 75252 Paris Cedex 05, France.}

\maketitle
\begin{abstract}
\noindent Exponentiating the hypergeometric series  $\:_0F_L(1,1,\ldots,1;z)$,
$L=0,1,2,\ldots,$ furnishes a recursion relation for the members of
integer sequences $b_L(n)$, $n=0,1,2,\ldots$. For $L>0$, the $b_L(n)$'s are
certain generalizations of conventional Bell numbers, $b_0(n)$.
The corresponding associated Stirling numbers of the second kind are
also generated and investigated. For $L=1$ one can give a
combinatorial interpretation of the numbers $b_1(n)$, and of some
Stirling numbers associated with them. We also consider  the $L\geq1$
analogues of Bell numbers for restricted partitions.
\end{abstract}

\bigskip

The conventional Bell numbers \cite{yablonsky}  $b_0(n)$, $n=0,1,2,\ldots$, have a well known exponential generating function
\begin{equation}
B_0(z)\equiv e^{\textstyle{(e^z-1)}}=\sum_{n=0}^{\infty}b_0(n)\frac{z^n}{n!},\label{B}
\end{equation}
which can be derived by interpreting $b_0(n)$ as the number of partitions of a set of $n$ distinct elements.
In this note we obtain recursion relations for related sequences of
positive integers, called $b_L(n)$, $L=0,1,2,\ldots,$ obtained by exponentiating the hypergeometric
series $\:_0F_L(1,1,\ldots,1;z)$ defined by \cite{andrews}:
\begin{equation}
\:_0F_L(\underbrace{1,1,\ldots,1}_{L};z)=\sum_{n=0}^{\infty}\frac{z^n}{(n!)^{L+1}},\label{defF}
\end{equation}
(for which we shall use throughout the short notation $\:_0F_L(z)$) and
which includes the special cases $\:_0F_0(z)\equiv e^z$ and $\:_0F_1(z)\equiv
I_0(2\sqrt{z})$, where $I_0(x)$ is the modified Bessel function of the
first kind. For $L>1$, the functions 
$\:_0F_L(z)$ are related to the so-called hyper-Bessel functions
\cite{marichev}, \cite{kiryakova}, \cite{paris},  which have recently
found application in quantum mechanics \cite{witte},
\cite{klauder}. Thus, we are interested in $b_L(n)$ given by
\begin{equation}
e^{[\:_0F_L(z)-1]}=\sum_{n=0}^{\infty}b_L(n)\frac{z^n}{(n!)^{L+1}},\label{geneF}   
\end{equation}
thereby defining a \underline{hypergeometric} generating function for
the numbers $b_L(n)$. 
From eq.(\ref{geneF}) it follows formally that 
\begin{equation}
b_L(n)=(n!)^L\cdot\frac{d^n}{dz^n}\left.\left(e^{[\:_0F_L(z)-1]}\right)\right|_{z=0}.\label{deriv}
\end{equation}
For $L=0$ the r.h.s of eq.(\ref{deriv}) can be evaluated in closed form:
\begin{equation}
b_0(n)=\frac{1}{e}\sum_{k=0}^{\infty}\frac{k^n}{k!}=\left\{\frac{1}{e^z}\left[\left(z\frac{d}{dz}\right)^ne^z\right]\right\}_{z=1}.\label{dob}
\end{equation}
The first equality in (\ref{dob}) is the celebrated Dobi\'{n}ski formula
\cite{yablonsky}, \cite{comtet}, \cite{wilf}. The second equality in
eq.(\ref{dob}) follows from observing that for a power series
$R(z)=\sum_{k=0}^{\infty}A_kz^k$
\begin{equation}
\left(z\frac{d}{dz}\right)^nR(z)=\sum_{k=0}^{\infty}A_k\:k^n\:z^k\label{opderiv}
\end{equation}
holds, and by applying eq.(\ref{opderiv}) to the exponential series $(A_k=(k!)^{-1})$.

The reason for  including the divisors $(n!)^{L+1}$ rather than $n!$ as in the usual exponential generating function, arises from the fact that only through eq.(\ref{geneF}) are the numbers $b_L(n)$  actually integers.
This can be seen from general formulas for exponentiation of a power
series \cite{comtet}, which employ the (exponential) Bell polynomials,
complicated and rather unwieldy objects. It cannot however be
considered as a proof that the $b_L(n)$ are integers.  At
this stage we shall use  eq.(\ref{geneF}) with $b_L(n)$ real and apply to it an
efficient method, exposed in \cite{wilf}, which will yield the recursion
relation for the $b_L(n)$. (For the proof that the $b_L(n)$ are integers, see
below eq.(\ref{recur})).   To this end we first obtain a result for the
multiplication of two power-series of the type (\ref{geneF}). Suppose
that we have to multiply $f(x)=\sum_{n=0}^{\infty}a_L(n)\frac{x^n}{(n!)^{L+1}}$ and
$g(x)=\sum_{n=0}^{\infty}c_L(n)\frac{x^n}{(n!)^{L+1}}$. We get
$f(x)\cdot g(x)=\sum_{n=0}^{\infty}d_L(n)\frac{x^n}{(n!)^{L+1}}$, where
\begin{equation}
d_L(n)=(n!)^{L+1}\sum_{r+s=n}^{\infty}\frac{a_L(r)c_L(s)}{(r!)^{L+1}(s!)^{L+1}}=\sum_{r=0}^{n}\left(\scriptsize{\begin{array}{c} n \\ r \end{array}}\right)^{L+1}a_L(r)\:c_L(n-r).\label{dL}
\end{equation}
Substitute eq.(\ref{defF}) into eq.(\ref{geneF}) and take the logarithm of both
sides of eq.(\ref{geneF}):
\begin{equation}
\sum_{n=1}^{\infty}\frac{z^n}{(n!)^{L+1}}= \ln\left(\sum_{n=0}^{\infty}b_L(n)\frac{z^n}{(n!)^{L+1}}\right).\label{log}
\end{equation}
Now differentiate both sides of eq.(\ref{log}) and multiply  by $z$. It produces
\begin{equation}
\left(\sum_{n=0}^{\infty}b_L(n)\frac{z^n}{(n!)^{L+1}}\right)\left(\sum_{n=0}^{\infty}n\:\frac{z^n}{(n!)^{L+1}}\right)=\sum_{n=0}^{\infty}n\:b_L(n)\frac{z^n}{(n!)^{L+1}},
\end{equation}
which with eq.(\ref{dL}) yields the desired recurrence relation
\begin{eqnarray}
b_L(n+1) & = &
\frac{1}{n+1}\sum_{k=0}^{n}\scriptsize{\left(\!\!\begin{array}{c} n+1 \\ k
\end{array}\!\!\right)}^{L+1}(n+1-k)\:b_L(k),\hspace{1cm}n=0,1,\ldots\label{recurb}\\
& = & \sum_{k=0}^{n}\left(\scriptsize{\begin{array}{c} n \\ k\end{array}}\right)\scriptsize{\left(\begin{array}{c} n+1 \\ k\end{array}\right)}^{L}\:b_L(k),\label{recur}\\
b_L(0) & = & 1.\label{bLinit}
\end{eqnarray}
Since eq.(\ref{recur}) involves only positive integers, it follows
that the $b_L(n)$ are indeed positive integers.
For  $L=0$ one gets the known recurrence relation for the Bell numbers \cite{wilf}:
\begin{equation}
b_0(n+1) =\sum_{k=0}^{n}\scriptsize{\left(\begin{array}{c} n \\ k \end{array}\right)}b_0(k).
\end{equation}
We have used eq.(\ref{recur}) to calculate some of the $b_L(n)$'s,
listed in Table I, for $L=0,1,\ldots,6$. Eq.(\ref{recur}), for $n$ fixed, gives 
closed form expressions for the $b_L(n)$ directly as a function of $L$ (columns
in Table I):  $b_L(2)=1+2^L$,  $b_L(3)=1+3\cdot3^L+(3!)^L$,
$b_L(4)=1+4\cdot4^L+3\cdot6^L+6\cdot12^L+(4!)^L$, etc.

The sets of $b_L(n)$ have been checked against the
most complete source of integer sequences available \cite{sloane}. Apart from the case
$L=0$ (conventional Bell numbers) only the first non-trivial sequence
$L=1$ is listed: it turns out that this sequence $b_1(n)$, listed under the heading
\underline{A023998} in \cite{sloane}, can be given a combinatorial
interpretation as the number of block permutations on a set of $n$ objects, which are
uniform, i.e. corresponding blocks have the same size
\cite{fitzgerald}.

Eq.(\ref{B}) can be generalized by including an additional
variable $x$, which will result in ``smearing out'' the conventional
Bell numbers $b_0(n)$ with a set of integers $S_0(n,k)$, such that for $k>n$,
$S_0(n,k)=0$, and $S_0(0,0)=1$, $S_0(n,0)=0$. In particular,
\begin{equation}
B_0(z,x)\equiv e^{\textstyle x(e^z-1)}=\sum_{n=0}^{\infty}\left[\sum_{k=1}^{n}S_0(n,k)\:x^k\right]\frac{z^n}{n!},\label{B0}
\end{equation}
which leads to the (exponential) generating function of $S_0(n,l)$,
the conventional Stirling numbers of the second kind, (see \cite{yablonsky}, \cite{comtet}), in the form
\begin{equation}
\frac{(e^z-1)^l}{l!}=\sum_{n=l}^{\infty}\frac{S_0(n,l)}{n!}z^n,
\end{equation}
and defines the so-called exponential or Touchard polynomials $l_n^{(0)}(x)$ as
\begin{equation}
l_n^{(0)}(x)=\sum_{k=1}^{n}S_0(n,k)x^k.
\end{equation}
They satisfy
\begin{equation}
l_n^{(0)}(1)=b_0(n),
\end{equation}
justifying  the term ``smearing out'' used above. 

The appearance of integers in eq.(\ref{geneF}) suggests a  natural
extension with an additional variable $x$:
\begin{equation}
B_L(z,x)\equiv e^{x[\:_0F_L(z)-1]}=\sum_{n=0}^{\infty}\left[\sum_{k=1}^{n}S_L(n,k)\:x^k\right]\frac{z^n}{(n!)^{L+1}},\label{expF}
\end{equation}
where we  include the right divisors $(n!)^{L+1}$
in the r.h.s of (\ref{expF}).

This in turn defines  ``hypergeometric''polynomials
of type $L$ and order $n$ through
\begin{equation}
l_n^{(L)}(x)=\sum_{k=1}^{n}S_L(n,k)x^k,\label{poly}
\end{equation}
which satisfy
\begin{equation}
l_n^{(L)}(1)=b_L(n),
\end{equation}
with the $b_L(n)$ of eq.(\ref{recurb}). Thus, the polynomials of
eq.(\ref{poly}) "smear out" the $b_L(n)$ with the generalized Stirling
numbers of the second kind, of type $L$, denoted by $S_L(n,k)$ (with
$S_L(n,k)=0$, if $k>n$, $S_L(n,0)=0$ if $n>0$ and $S_L(0,0)=1$), which have, from eq.(\ref{expF}) the ``hypergeometric''generating function
\begin{equation}
\frac{(\:_0F_L(z)-1)^l}{l!}=\sum_{n=l}^{\infty}\frac{S_L(n,l)}{(n!)^{L+1}}\:z^n,\hspace{1cm}L=0,1,2,\ldots .\label{geneSL}
\end{equation}

Eq.(\ref{geneSL}) can be used to derive a recursion relation for
the numbers $S_L(n,k)$, in the same manner as eq.(\ref{geneF})
yielded eq.(\ref{bLinit}). Thus we take the logarithm of both sides
of eq.(\ref{geneSL}), differentiate with respect to $z$, multiply by $z$
and obtain:
\begin{equation}
\left(\sum_{n=0}^{\infty}\frac{S_L(n,l-1)}{(n!)^{L+1}}\:z^n\right)\left(\sum_{n=0}^{\infty}\frac{n}{(n!)^{L+1}}\:z^n\right)=\sum_{n=0}^{\infty}\frac{n\:S_L(n,l)}{(n!)^{L+1}}\:z^n,\label{eqsnl}
\end{equation}
which, with the help of eq.(\ref{dL}),  produces the required recursion
relation
\begin{eqnarray}
S_L(n+1,l) & = & \sum_{k=l-1}^{n}\left(\scriptsize{\begin{array}{c} n \\
k\end{array}}\right)\scriptsize{\left(\begin{array}{c} n+1 \\
k\end{array}\right)}^{L}\:S_L(k,l-1),\label{recurSL}\\
& &  ,\hspace{-1cm}S_L(0,0)=1,\hspace{1cm}S_L(n,0)=0,\label{initrecurSL}
\end{eqnarray}
which for $L=0$ is the recursion relation for the conventional
Stirling numbers of the second kind \cite{yablonsky}, \cite{comtet}, and
in eq.(\ref{recurSL}) the appropriate summation range has been inserted.
Since the recursions of eq.(\ref{recurSL}) and eq.(\ref{initrecurSL}) involve
only integers we conclude that $S_L(n,l)$ are positive integers.

We have calculated some of the numbers $S_L(n,l)$ using
eq.(\ref{geneSL}) and have listed them in Tables II and III,
for $L=1$ and $L=2$ respectively. Observe that $S_1(n,2)=\left(\scriptsize{\begin{array}{c} 2n+1 \\
n+1\end{array}}\right)-1$ and $S_L(n,n)=(n!)^L$, $L=1,2$. Also, by
fixing $n$ and $l$, the individual values of $S_L(n,l)$ have been
calculated as a function of $L$ with the help of eq.(\ref{recurSL}),
see Table IV, from which we observe
\begin{equation}
S_L(n,n)=(n!)^L,\hspace{1cm} L=1,2,\ldots.\label{Snn}
\end{equation}
which is the lowest diagonal in Table IV. We now demonstrate that the
repetitive use of eq.(\ref{recurSL}) permits one to establish
closed-form expressions for any supra-diagonal of order $p$, i.e. the sequence $S_L(n+p,n)$, for $p=1,2,3,\ldots$,
if one knows the expression for all $S_L(n+k,n)$ with $k<p$.
We shall illustrate it here for $p=1,2$. To this end fix $l=n$ on both
sides of eq.(\ref{recurSL}). It becomes, upon using eq.(\ref{Snn}),
and defining $\alpha_L(n)\equiv S_L(n+1,n)$, a linear recursion relation
\begin{equation}
\alpha_L(n)=\frac{n[(n+1)!]^L}{2^L}+(n+1)^L\alpha_L(n-1),\hspace{1cm}\alpha_L(0)=0,
\end{equation}
with the solution
\begin{eqnarray}
\alpha_L(n)=S_L(n+1,n)&=&\frac{n(n+1)}{2}\left[\frac{(n+1)!}{2}\right]^L\label{alpha}\\
&=&\left[\frac{(n+1)!}{2}\right]^L\:S_0(n+1,n),\label{s1}
\end{eqnarray}
which gives the second lowest diagonal in Table IV. Observe that for
any $L$, $S_L(n+1,n)$ is proportional to $S_0(n+1,n)=n(n+1)/2$. The
sequence $S_1(n+1,n)=1,\:9,\:72,\:600,\:5\:400,8\:564\:480,\:\ldots$ is of
particular interest: it represents the sum of inversion numbers of all
permutations on $n$ letters \cite{sloane}. For more information
about this and related sequences see the entry \underline{A001809} in \cite{sloane}.
The $S_L(n+1,n)$ for $L>1$ do not appear to have a simple combinatorial
interpretation.
A recurrence equation for $\beta_L(n)\equiv S_L(n+2,n)$ is obtained upon
substituting eq.(\ref{Snn}) and eq.(\ref{alpha}) into eq.(\ref{recurSL}):
\begin{equation}
\beta_L(n)=\frac{n(n+1)}{2!}\left[\frac{(n+2)!}{2!}\right]^L\left(\frac{n-1}{2^L}+\frac{1}{3^L}\right)+(n+2)^L\beta_L(n-1),\hspace{1cm}\beta_L(0)=0.
\end{equation}
It has the solution
\begin{equation}
S_L(n+2,n)=\frac{n(n+1)(n+2)}{3\cdot2^3}\left[\frac{(n+2)!}{2}\right]^L\left(\frac{3}{2^L}(n-1)+\frac{4}{3^L}\right)\label{SL}
\end{equation}
which is a closed form expression for the second lowest diagonal in
Table IV. Clearly, eq.(\ref{SL}) for $L=0$ gives the combinatorial
form for the series of conventional Stirling numbers
\begin{equation}
S_0(n+2,n)=\frac{n(n+1)(n+2)(3n+1)}{4!}.\label{s2}
\end{equation}
In a similar way we obtain
\begin{eqnarray}
S_L(n+3,n)&=&\frac{n(n+1)(n+2)(n+3)}{3\cdot2^4}\left[\frac{(n+3)!}{3}\right]^L\nonumber\\
& &\times\left(n^2\left(\frac{3}{8}\right)^L+n\left(\frac{1}{4^{L-1}}-\frac{3^{L+1}}{8^L}\right)+\frac{2+2\cdot3^L}{8^L}-\frac{1}{4^{L-1}}\right)
\end{eqnarray}
which for $L=0$ reduces to
\begin{equation}
S_0(n+3,n)=\frac{1}{48}n^2(n+1)^2(n+2)(n+3).\label{s3}
\end{equation}
Combined with the standard definition \cite{comtet}, \cite{wilf}
\begin{equation}
S_0(n,l)=\frac{(-1)^l}{l!}\sum_{k=1}^{l}(-1)^k\left(\scriptsize{\begin{array}{c} l \\
k\end{array}}\right)\:k^n.\label{A1}
\end{equation}
eqs.(\ref{s1}), (\ref{s2}) and (\ref{s3}) give compact expressions for
the summation form of $S_0(n+p,n)$. Further, from eq.(\ref{A1}), use
of eq.(\ref{opderiv}) gives the following generating formula
\begin{eqnarray}
S_0(n,l) & = &\frac{(-1)^l}{l!}\left[\left(z\frac{d}{dz}\right)^n\left(\sum_{k=1}^{l}(-1)^k\left(\scriptsize{\begin{array}{c}l \\ k\end{array}}\right)\:z^k\right)\right]_{z=1} \\
& = & \frac{(-1)^l}{l!}\left[\left(z\frac{d}{dz}\right)^n[(1-z)^l-1]\right]_{z=1},\hspace{1cm}n\geq l.
\end{eqnarray}

The formula (\ref{B}) can be generalized by putting restrictions on
the type of resulting partitions. The generating function for the
number of partitions of a set of $n$ distinct  elements without
singleton blocks $b_0(1,n)$ is \cite{comtet}, \cite{ehrenborg}, \cite{suter},
\begin{equation}
B_0(1,z)=e^{e^{z}-1-z}=\sum_{n=0}^{\infty}b_0(1,n)\frac{z^n}{n!},
\end{equation}
or more generally, without singleton, doubleton $\ldots$, $p-$blocks
$(p=0,1,\ldots)$ is \cite{suter}
\begin{equation}
B_0(p,z)=e^{e^{z}-\sum_{k=0}^{p}\frac{z^k}{k!}}=\sum_{n=0}^{\infty}b_0(p,n)\frac{z^n}{n!},
\end{equation}
with the corresponding associated Stirling numbers defined by analogy with
eq.(\ref{B0}) and eq.(\ref{eqsnl}).
The numbers $b_0(1,n)$, $b_0(2,n)$, $b_0(3,n)$, $b_0(4,n)$ can be read
off from the sequences \underline{A000296},  \underline{A006505},
\underline{A057837} and  \underline{A057814} in \cite{sloane},
respectively. For more properties of these numbers see \cite{bernstein}.

We carry over this type of extension to eq.(\ref{geneF}) and define $b_L(p,n)$ through
\begin{equation}
B_L(p,z)\equiv e^{\:_0F_L(z)-\sum_{k=0}^{p}\frac{z^k}{(k!)^{L+1}}}=\sum_{n=0}^{\infty}b_L(p,n)\frac{z^n}{(n!)^{L+1}},\label{BLpz}
\end{equation}
where $b_L(0,n)=b_L(n)$ from eq.(\ref{geneF}). (We know of no
combinatorial meaning of $b_L(p,n)$ for $L\geq1$, $p>0$).
The $b_L(p,n)$ satisfy the following recursion relations:
\begin{eqnarray}
b_L(p,n)& = & \sum_{k=0}^{n-p}\left(\scriptsize{\begin{array}{c} n \\ k\end{array}}\right)\scriptsize{\left(\begin{array}{c} n+1 \\ k\end{array}\right)}^{L}\:b_L(p,k),\label{defblpn}\\
b_L(p,0)& = & 1,\\
b_L(p,1)& = & b_L(p,2)=\cdots=b_L(p,p)=0,\\
b_L(p,p+1)& = & 1.
\end{eqnarray}
That the $b_L(p,n)$ are integers follows from eq.(\ref{defblpn}).
 Through eq.(\ref{BLpz}) additional families of integer Stirling-like
numbers $S_{L,p}(n,k)$ can be readily defined and investigated. 

The numbers $b_0(p,n)$ are collected in Table V, and Tables VI and
VII contain the lowest values of $b_1(p,n)$ and $b_2(p,n)$, respectively.

Formula (\ref{B}) can be used to express $e$ in terms of $b_0(n)$
in various ways. Two such lowest order (in differentiation) forms are
\begin{eqnarray}
e & = & 1+\ln\left( \sum_{n=0}^{\infty}\frac{b_0(n)}{n!}\right)=\label{lnb0} \\
& = &\ln\left( \sum_{n=0}^{\infty}\frac{b_0(n+1)}{n!}\right).\label{lnb01}
\end{eqnarray}   
In the very same way, eq.(\ref{geneF}) can be used to express the
values of $\:_0F_L(z)$ and its derivatives at $z=1$ in terms of certain
series of $b_L(n)$'s. For $L=1$, the analogues of eq.(\ref{lnb0}) and
eq.(\ref{lnb01}) are 
\begin{eqnarray}
I_0(2) & = & 1+\ln\left( \sum_{n=0}^{\infty}\frac{b_1(n)}{(n!)^2}\right), \\
I_0(2)+\ln(I_1(2)) & = & 1+\ln\left(\sum_{n=0}^{\infty}\frac{b_1(n+1)}{(n+1)(n!)^2}\right)
\end{eqnarray}
and for $L=2$ the corresponding formulas are
\begin{eqnarray}
\:_0F_2(1,1;1) & = & 1+\ln\left( \sum_{n=0}^{\infty}\frac{b_2(n)}{(n!)^3}\right), \\
\:_0F_2(1,1;1)+\ln\left(\:_0F_2(2,2;1)\right) & = & 1+\ln\left(\sum_{n=0}^{\infty}\frac{b_2(n+1)}{(n+1)^2(n!)^3}\right).
\end{eqnarray}

By fixing  $z_0$ at  values other than $z_0=1$, one can link the numerical
values of certain combinations of $\:_0F_L(1,1,\ldots;z_0)$ ,
$\:_0F_L(2,2,\ldots;z_0)$,\ldots and their logarithms, with other series
containing the $b_L(n)$'s.

The above considerations can be extended to the exponentiation of the more
general hypergeometric functions of type $\:_0F_L(k_1,k_2,\ldots,k_L;z)$
where $k_1,k_2,\ldots,k_L$ are positive integers. We conjecture that for every set of $k_n$'s
a different set of integers will be generated through an appropriate
adaptation of eq.(\ref{geneF}). We quote one simple example of such a
series. For
\begin{equation}
\:_0F_2(1,2;z)=\sum_{n=0}^{\infty}\frac{z^n}{(n+1)(n!)^3}
\end{equation}
eq.(\ref{geneF}) extends to
\begin{equation}
e^{[\:_0F_2(1,2;z)-1]}=\sum_{n=0}^{\infty}f_2(n)\frac{z^n}{(n+1)(n!)^3}
\end{equation}
where 
\begin{equation}
f_2(n)=(n+1)(n!)^2\left[\frac{d^n}{dz^n}e^{[\:_0F_2(1,2;z)-1]}\right]_{z=0}
\end{equation}
turn out to be integers:
$f_2(n)$, $n=0,1,\ldots,8$ are: 1, 1, 4, 37, 641, 18 276, 789 377, 48 681
011, etc.
The analogue of equations (\ref{recurSL}) and (\ref{lnb0}) is:
\begin{equation}
\:_0F_2(1,2;1)=1+\ln\left(\sum_{n=0}^{\infty}\frac{f_2(n)}{(n+1)(n!)^3}\right).
\end{equation}

\bigskip

\acknowledgements
We thank L. Haddad for interesting discussions. We have used ${\rm
Maple}^{\copyright}$ to calculate most of the numbers discussed above.

\begin{table}
\caption{Table of $b_L(n)$; $L,n=0,1,\ldots,6.$}
\begin{tabular}{cccccccc}
$L$ & $b_L(0)$ & $b_L(1)$ & $b_L(2)$ & $b_L(3)$ & $b_L(4)$ & $b_L(5)$ & $b_L(6)$ \\ \tableline
0 &\hspace{0.3cm} 1 & 1 & 2 & 5 & 15 & 52 & 203 \\
1 &\hspace{0.3cm} 1 & 1 & 3 & 16 & 131 & 1 496 & 22 482 \\
2 &\hspace{0.3cm} 1 & 1 & 5 & 64 & 1 613 & 69 026 & 4 566 992 \\
3 &\hspace{0.3cm} 1 & 1 & 9 & 298 & 25 097 & 4 383 626 & 1 394 519 922\\
4 &\hspace{0.3cm} 1 & 1 & 17 & 1 540 & 461 105 & 350 813 126 & 573 843
627 152\\ 
5 &\hspace{0.3cm} 1 & 1 & 33 & 8 506 & 9 483 041 & 33 056 715 626 &
293 327 384 637 282\\
6 &\hspace{0.3cm} 1 & 1 & 65 & 48 844 & 209 175 233 & 3 464 129 078
126 & 173 566 857 025 139 312\\
\end{tabular}
\end{table}

\begin{table}
\caption{Table of $S_L(n,l)$;  $L=1$,  $l,n=1,2,\ldots,8.$}
\begin{tabular}{ccccccccc}
$l$ & $S_1(1,l)$ & $S_1(2,l)$ & $S_1(3,l)$ & $S_1(4,l)$ & $S_1(5,l)$ &
$S_1(6,l)$ & $S_1(7,l)$ & $S_1(8,l)$ \\ \tableline
1 &\hspace{0.3cm} 1 & 1 & 1 & 1 & 1 & 1 & 1 & 1\\
2 &\hspace{0.3cm}   & 2 & 9 & 34 & 125 & 461 & 1 715 & 6 434 \\
3 &\hspace{0.3cm}   &   & 6 & 72 & 650 & 5 400 & 43 757 & 353 192 \\
4 &\hspace{0.3cm}   &   &   & 24 & 600 & 10 500 & 161 700 & 2 361 016 \\
5 &\hspace{0.3cm}   &   &   &    & 120 & 5 400 & 161 700 & 4 116 000  \\ 
6 &\hspace{0.3cm}   &   &   &    &     & 720  & 52 920 & 2 493 120 \\
7 &\hspace{0.3cm}   &   &   &    &     &      & 5 040 & 564 480 \\
8 &\hspace{0.3cm}   &   &   &    &     &      &     & 40 320 \\
\end{tabular}
\end{table}

\begin{table}
\caption{Table of $S_L(n,l)$;  $L=2$,  $l,n=1,2,\ldots,8.$}
\begin{tabular}{ccccccccc}
$l$ & $S_2(1,l)$ & $S_2(2,l)$ & $S_2(3,l)$ & $S_2(4,l)$ & $S_2(5,l)$ &
$S_2(6,l)$ & $S_2(7,l)$ & $S_2(8,l)$ \\ \tableline
1 &\hspace{0.3cm} 1 & 1 & 1 & 1 & 1 & 1 & 1 & 1\\
2 &\hspace{0.3cm}   & 4 & 27 & 172 & 1 125 & 7 591 & 52 479 & 369 580 \\
3 &\hspace{0.3cm}   &   & 36 & 864 & 17 500 & 351 000 & 7 197 169 & 151 633 440 \\
4 &\hspace{0.3cm}   &   &   & 576 & 36 000 & 1 746 000 & 80 262 000 &
3 691 514 176\\
5 &\hspace{0.3cm}   &   &   &    & 14 400 & 1 944 000 & 191 394 000 & 17 188 416 000\\ 
6 &\hspace{0.3cm}   &   &   &    &     & 518 400  & 133 358 400 & 23
866 214 400 \\
7 &\hspace{0.3cm}   &   &   &    &     &    & 25 401 600   & 11 379 916 800\\
8 &\hspace{0.3cm}   &   &   &    &     &    &            & 1 625 702 400\\
\end{tabular}
\end{table}

\begin{table}
\caption{Table of $S_L(n,l)$;  $l,n=1,2,\ldots,6.$}
\begin{tabular}{ccccccc}
$l$ & $S_L(1,l)$ & $S_L(2,l)$ & $S_L(3,l)$ & $S_L(4,l)$ & $S_L(5,l)$ &
$S_L(6,l)$ \\ \tableline
1 &\hspace{0.3cm} 1 & 1 & 1 & 1 & 1 & 1\\
2 &\hspace{0.3cm}   & $(2!)^L$ & $3\cdot3^L$ & $4\cdot4^L+3\cdot6^L$ & $5\cdot5^L+10\cdot10^L$ & $6\cdot6^L+15\cdot15^L+10\cdot20^L$\\
3 &\hspace{0.3cm}   &   & $(3!)^L$ & $6\cdot12^L$ & $10\cdot20^L$+$15\cdot30^L$
& $15\cdot30^L+60\cdot60^L+15\cdot90^L$ \\
4 &\hspace{0.3cm}   &   &   & $(4!)^L$ & $10\cdot60^L$ & $20\cdot120^L+45\cdot180^L$\\
5 &\hspace{0.3cm}   &   &   &    & $(5!)^L$ & $15\cdot360^L$ \\ 
6 &\hspace{0.3cm}   &   &   &    &  & $(6!)^L$ \\ 
\end{tabular}
\end{table}

\begin{table}
\caption{Table of $b_0(p,n)$; $p=0,1,2,3;\:\:n=0,\ldots,10.$}
\begin{tabular}{ccccc}
$n$ & $b_0(0,n)$ & $b_0(1,n)$ & $b_0(2,n)$ & $b_0(3,n)$ \\ \tableline
0 &\hspace{0.3cm} 1 & 1 & 1 & 1 \\
1 &\hspace{0.3cm} 1 & 0 & 0 & 0 \\
2 &\hspace{0.3cm} 2 & 1 & 0 & 0 \\
3 &\hspace{0.3cm} 5 & 1 & 1 & 0 \\
4 &\hspace{0.3cm} 15 & 4 & 1 & 1 \\ 
5 &\hspace{0.3cm} 52 & 11 & 1 & 1 \\
6 &\hspace{0.3cm} 203 & 41 & 11 & 1\\
7 &\hspace{0.3cm} 877 & 162 & 36 & 1 \\
8 &\hspace{0.3cm} 4 140 & 715 & 92 & 36 \\
9 &\hspace{0.3cm} 21 147 & 3 425 & 491 & 127 \\
10 &\hspace{0.3cm} 115 975 & 17 722 & 2 557 & 337 \\
\end{tabular}
\end{table}

\begin{table}
\caption{Table of $b_1(p,n)$; $p=0,1,2;\:\:n=0,\ldots,9.$}
\begin{tabular}{cccc}
$n$ & $b_1(0,n)$ & $b_1(1,n)$ & $b_1(2,n)$ \\ \tableline
0 &\hspace{0.3cm} 1 & 1 & 1 \\
1 &\hspace{0.3cm} 1 & 0 & 0 \\
2 &\hspace{0.3cm} 3 & 1 & 0 \\
3 &\hspace{0.3cm} 16 & 1 & 1\\
4 &\hspace{0.3cm} 131 & 19 & 1 \\ 
5 &\hspace{0.3cm} 1 496  & 101 & 1 \\
6 &\hspace{0.3cm} 22 482 & 1 776 &  201\\
7 &\hspace{0.3cm} 426 833 & 23 717 & 1 226 \\
8 &\hspace{0.3cm} 9 934 563 & 515 971 & 5 587 \\
9 &\hspace{0.3cm} 277 006 192 & 11 893 597 & 493 333 \\
\end{tabular}
\end{table}

\begin{table}
\caption{Table of $b_2(p,n)$; $p=0,1,2;\:\:n=0,\ldots,8.$}
\begin{tabular}{cccc}
$n$ & $b_2(0,n)$ & $b_2(1,n)$ & $b_2(2,n)$ \\ \tableline
0 &\hspace{0.3cm} 1 & 1 & 1 \\
1 &\hspace{0.3cm} 1 & 0 & 0 \\
2 &\hspace{0.3cm} 5 & 1 & 0 \\
3 &\hspace{0.3cm} 64 & 1 & 1\\
4 &\hspace{0.3cm} 1 613 & 109 & 1 \\ 
5 &\hspace{0.3cm} 69 026  & 1 001 & 1 \\
6 &\hspace{0.3cm} 4 566 992 & 128 876 &  4 001\\
7 &\hspace{0.3cm} 437 665 649 & 4 682 637 & 42 876 \\
8 &\hspace{0.3cm} 57 903 766 800 & 792 013 069 & 347 117 \\
\end{tabular}
\end{table}


\begin{thebibliography}{99}      
\bibitem{yablonsky}S.V. Yablonsky, `` Introduction to Discrete
Mathematics'', (Mir Publishers, Moscow, 1989)
\bibitem{andrews} G.E. Andrews, R. Askey and R. Roy, ``Special
Functions'', Encyclopedia of Mathematics and its Applications, vol.71,
(Cambridge University Press,  1999)
\bibitem{marichev} O.I. Marichev, \emph{Handbook of Integral
Transforms of Higher Transcendental Functions, Theory and Algorithmic
Tables}, (Ellis Horwood Ltd, Chichester, 1983), Chap. 6
\bibitem{kiryakova}V.S. Kiryakova and B.Al-Saqabi, ``Explicit solutions to
hyper-Bessel integral equations of second kind'', Comput. and
Math. with Appl. \textbf{37}, 75 (1999)
\bibitem{paris}R.B. Paris and A.D. Wood, ``Results old and new on the
hyper-Bessel equation'', Proc.Roy.Soc. Edinb. \textbf{106A}, 259 (1987)
\bibitem{witte}N.S. Witte, ``Exact solution for the reflection and
diffraction of atomic de Broglie waves by a traveling evanescent laser
wave'', J.Phys.A\textbf{31}, 807 (1998)
\bibitem{klauder}J.R. Klauder, K.A. Penson and J.-M. Sixdeniers,
``Constructing coherent states through solutions of Stieltjes and
Hausdorff moment problems'', (The Physical Review A, in press (2001))
\bibitem{comtet}L. Comtet, "Advanced Combinatorics", (D. Reidel,
Boston, 1984)
\bibitem{wilf} H.S. Wilf, ``Generatingfunctionology'', $2^{\rm nd}$
ed., (Academic Press, New York, 1994)
\bibitem{sloane}N.J.A. Sloane, On-Line Encyclopedia of Integer
Sequences, published electronically at:
http://www.research.att.com/$\sim$/njas/sequences/
\bibitem{bernstein}M. Bernstein and N.J.A. Sloane, "Some canonical
sequences of integers'', Linear Algebra Appl., \textbf{226/228}, 57 (1995)
\bibitem{fitzgerald}D.G. Fitzgerald and J. Leech, ``Dual symmetric
inverse monoids and representation theory'', J.Austr.Math.Soc., Series
A\textbf{64}, 345 (1998)
\bibitem{delerue}P. Delerue, ``Sur le calcul symbolique \`{a} $n$
variables et fonctions hyperbess\'{e}liennes II'',
Ann.Soc.Sci. Brux. \textbf{67}, 229 (1953)
\bibitem{ehrenborg}R. Ehrenborg, "The Hankel Determinant of Exponential
Polynomials", Am.Math.Monthly,  \textbf{207}, 557 (2000)
\bibitem{suter} R. Suter, ``Two Analogues of a Classical Sequence",
J.Integ.Seq.  \textbf{3}, Article 00.1.8 (2000), available
electronically through \cite{sloane}
\end{thebibliography}
\end{document}